\documentclass[12pt,oneside]{amsart}
\usepackage[margin=1in]{geometry}
\usepackage{graphicx, amssymb, tikz, latexsym, epsfig, enumitem, multicol, wasysym, amsmath, amsthm, amscd, amssymb, latexsym, ulem, soul, multirow, rotating, pdflscape, float, placeins, dsfont, comment}
\usepackage[colorlinks,citecolor=red,pagebackref,hypertexnames=false]{hyperref}
\hypersetup{%
  colorlinks = true,
  citecolor=red,
  linkcolor  = red
}
\setlist{itemsep=.06125in}
 
\usepackage{graphicx}
\restylefloat{table}
\numberwithin{equation}{section}
\theoremstyle{plain}
\newtheorem{theorem}{Theorem}[section]
\newtheorem{lemma}[theorem]{Lemma}

\theoremstyle{definition}
\newtheorem{definition}[theorem]{Definition}

\theoremstyle{remark}
\newtheorem{remark}[theorem]{Remark}

\newtheorem{case[theorem]}{Case}

\date{\today}      

\thanks{A.~I. was supported in part by the National Science Foundation under NSF DMS - 2154232. A.~M. was supported in part by AMS-Simons Research Enhancement Grant and  by the National Science Foundation under NSF DMS-2453769.}

\title{Fourier minimization and time series imputation}

\author{W. Burstein$^1$, A. Iosevich$^1$, A. Mayeli$^2$, and H. Nathan$^{1,*}$}

\begin{document}

\maketitle

\vspace{-1.5\baselineskip}

\begin{center}
    $^1$Department of Mathematics, University of Rochester, Rochester, NY, USA\\
    \textit{Emails:} wburste2@ur.rochester.edu; iosevich@gmail.com\\
    \vspace{0.15cm}
    $^2$Department of Mathematics, CUNY Graduate Center, New York, NY, USA\\
    \textit{Email:} amayeli@gc.cuny.edu\\
    \vspace{0.15cm}
    $^*$Corresponding author: \textit{Email:} hnathan3@ur.rochester.edu\\
\end{center}

\begin{abstract} One of the most common procedures in modern data analytics is filling in missing values in times series. For a variety of reasons, the data provided by clients to obtain a forecast, or other forms of data analysis, may have missing values, and those values need to be filled in before the data set can be properly analyzed. Many freely available forecasting software packages, such as the sktime library, have built-in mechanisms for filling in missing values. The purpose of this paper is to adapt the classical $L^1$ minimization method for signal recovery to the filling of missing values in time series. The theoretical justifications of these methods leverage results by Bourgain (\cite{Bourgain89}), Talagrand (\cite{Talagrand98}), the second and the third listed authors (\cite{IM24}), and the result by the second listed author, Kashin, Limonova and the third listed author (\cite{IKLM24}). Brief numerical tests for these algorithms are given but more extensive tests will be discussed in a companion paper.
\end{abstract}

\tableofcontents

\section{Introduction and main results}
\label{section:introduction}

The main object we study in this paper is the time series $f: {\mathbb Z}_N \to {\mathbb C}$, where ${\mathbb Z}_N=\{0,1,2, \dots, N-1\}$ denotes integers modulo $N$ which serves as a mathematical structure to encode signals of length $N$. Suppose that the set 

$$ \{f(x): x \in M\}$$

is missing, where $M \subset {\mathbb Z}_N$. The main thrust of this paper is to provide algorithms for recovering these values and to discuss the theoretical and applied justifications and limitations of these methods. The algorithms addressed in this paper try to find a function $g:\mathbb{Z}_N \to {\mathbb C}$, which agrees with $f$ on $M^c$ and, otherwise, minimizes the $L^1$-norm of the Fourier transform of $g$.

\vskip.125in 

This paper consists of theoretical results; readers interested in practical considerations can find them in the companion paper. The main results of this paper are Theorem \ref{theorem:looseLogan} and Theorem \ref{theorem:quantitativerandomnoisyLogan}. Theorem \ref{theorem:quantitativerandomnoisyLogan} describes the accuracy of imputation using Logan's $L^1$ minimization method when the set of missing values is randomly chosen and the data contains a certain amount of random noise. Theorem \ref{theorem:looseLogan} describes the situation where we make no assumptions on the set of missing values. 

\vskip.125in 

\subsection{Notation and Conventions}
\label{subsection:notation}

In this sub-section, we describe the notation and conventions used throughout this paper.

Our main tool is the classical discrete Fourier transform (see e.g. \cite{Smith07}). For a fixed $N$, let $\chi:\mathbb{Z}_N \to {\mathbb C}$ denote the standard character on ${\mathbb Z}_N$, i.e $\chi(s) := \exp\left( \frac{2 \pi i s}{N}\right)$. 

Given $f:{\mathbb Z}_N \to {\mathbb C}$, we define the discrete Fourier transform $\mod N$, $\widehat{f}:\mathbb{Z}_N \to \mathbb{C}$, by the formula 

$$\widehat{f}(\omega) := N^{-\frac{1}{2}} \sum_{x \in {\mathbb Z}_N} f(x) \chi(-\omega \cdot x).$$

The original function can be recovered from its Fourier transform. If $\widehat{f}$ is the Fourier transform of $f$ then for $x \in {\mathbb Z}_N$

$$f(x) = N^{-\frac{1}{2}} \sum_{\omega \in {\mathbb Z}_N} \widehat{f}(\omega) \chi(x \cdot \omega).$$

We will make frequent use of the fact that for any $s \in {\mathbb Z}_N$, $|\chi(s)| = 1$.

When discussing the $L_p$ norm of a function $f:X \to \mathbb{C}$ (for finite $X$) we define the un-normalized norm as

$$\| f \|_p = \|f\|_{L^p(X)} = \left(\sum_{x \in X} |f(x)|^p\right)^{1/p}$$

and the normalized norm as

$$\|f\|_{L^p(\mu)} = \left(\frac{1}{|X|}\sum_{x \in \mathbb{X}} |f(x)|^p\right)^{1/p}.$$

When we state that a function $f$ (or it's Fourier transform, $\widehat{f}$) is transmitted, we simply mean that some agent receives $f$ (or $\widehat{f}$) with some values possibly missing. For example, we an agent might receive the following data for a function $f:Z_4 \to \mathbb{C}$

\begin{center}
    \begin{tabular}{|c|c|c|c|c|}
        \hline
        $x$ & $0$ & $1$ & $2$ & $3$ \\
        \hline
        $f(x)$ & $1.1$ & $0.2$ & & $3.78$ \\
        \hline
    \end{tabular}
\end{center}

and would therefore know $f(0)$, $f(1)$, and $f(3)$ but not $f(2)$. We would say that ``$f$ was transmitted but $f(2)$ was unobserved.''

Throughout this paper we define $o_N(1)$ to be a quantity that goes to zero as $N \to \infty$. Wherever used, we assume that all other relevant parameters that are indepenent of $N$ are held constant as $N \to \infty.$

\vskip.125in 

\subsection{$L^1$-minimization algorithm and Logan's phenomenon} 
\label{subsection: L1 minimization and logan} Our algorithms are generalizations of the following result from Donoho and Stark (\cite{DS89}) as a consequence of the celebrated Logan's phenomenon.

\vskip.125in 

\begin{theorem}[\cite{DS89}]
\label{theorem:DSLogan89Orig} 
 Let $f: {\mathbb Z}_N \to {\mathbb C}$ be supported in $E \subset {\mathbb Z}_N$. Suppose that $\widehat{f}$ is transmitted but the frequencies ${\{\widehat{f}(m)\}}_{m \in S}$ are unobserved, where $S \subset {\mathbb Z}_N$, with $|E| \cdot |S|<\frac{N}{2}$. Then $f$ can be recovered exactly and uniquely. Moreover,

\begin{equation} 
\label{equation: simpleL1recoveryOrig}
f = argmin_g {||g||}_{L^1({\mathbb Z}_N)} \ \text{with the constraint} \ \widehat{f}(m)=\widehat{g}(m), m \notin S.
\end{equation} 
\end{theorem} 

The proof is remarkably illuminating and essentially proceeds by taking $h = f - g$ where $g$ is the solution of (\ref{equation: simpleL1recoveryOrig}). The triangle inequality gives us that

$${||g||}_{L^1({\mathbb Z}_N)} = {||f-h||}_{L^1(E)}+{||h||}_{L^1(E^c)} \geq {||f||}_{L^1({\mathbb Z}_N)}+\left({||h||}_{L^1(E^c)}-{||h||}_{L^1(E)} \right).$$

Using this, one can show that if $g \not\equiv f$ then, under the assumptions in the Theorem,

$${||h||}_{L^1(E^c)}-{||h||}_{L^1(E)} > 0 \Rightarrow \| g \|_{L^1(\mathbb{Z}_N)} > \| f \|_{L^1(\mathbb{Z}_N)}$$

which is a contradiction because $g$ has minimal $L^1$-norm. Full details can be found in \cite{DS89} and the Appendix (section \ref{section:appendix}) of this paper.

From theorem \ref{theorem:DSLogan89Orig} and the formula for Fourier inversion, we get the following theorem.

\begin{theorem} \label{theorem:DSLogan89} Let $f: {\mathbb Z}_N \to {\mathbb C}$ such that $\widehat{f}$ is supported in $E \subset {\mathbb Z}_N$. Suppose that $f$ is transmitted but the values ${\{f(x)\}}_{x \in M}$ are unobserved, where $M \subset {\mathbb Z}_N$, with $|E| \cdot |M|<\frac{N}{2}$. Then f can be recovered exactly and uniquely. Moreover,

\begin{equation} 
\label{equation: simpleL1recovery}
f = argmin_g {||\widehat{g}||}_{L^1({\mathbb Z}_N)} \ \text{with the constraint} \ f(x) = g(x), x \notin M.
\end{equation} 
\end{theorem}

\begin{remark}
Theorems \ref{theorem:DSLogan89Orig} and \ref{theorem:DSLogan89}, as well as the theorems below regarding $L^1$-minimization, are attractive for a variety of reasons. From a practical perspective, since the Fourier transform is a linear transformation, $L^1$-minimization can be implemented via convex optimization. As such, the algorithm is both efficient and relatively easy to implement.
\end{remark}

Our first result of this subsection is a quantitative version of Theorem \ref{theorem:DSLogan89Orig} and the main result in \cite{IKLM24}. 

\begin{definition} \label{def:eps_conc} Let $u: {\mathbb Z}_N \to {\mathbb C}$. We say that $u$ is $L^p$-concentrated on $A \subset {\mathbb Z}_N$ with norm $\leq \epsilon$ if 
$$ {||u||}_{L^p(A^c)} \leq \frac{\epsilon}{N} \cdot {||u||}_{L^p({\mathbb Z}_N)}.$$
\end{definition} 

\begin{theorem} \label{theorem:quantitativeLogan} Let $f: {\mathbb Z}_N \to {\mathbb C}$, and suppose that $f$ is transmitted but the values ${\{f(x)\}}_{x \in M}$ are unobserved. Suppose that $\widehat{f}$ is $L^1$-concentrated on $S \subset {\mathbb Z}_N$ with norm $\leq \epsilon$. Let $\delta=\frac{|M||S|}{N}$ and suppose that $\delta<\frac{1}{2}$. Let 
\begin{equation} \label{equation:simpleL1recoveryreversed} g = argmin_u {||\widehat{u}||}_{L^1({\mathbb Z}_N)} \ \ \text{with the constraint} \ f(x)=u(x) \ \text{for} \ x \notin M. \end{equation}  

Let $h=f-g$. Then if $h \not=0$, then 
\begin{equation} \label{equation:minerrorboundreversed} {||\widehat{h}||}_{L^1({\mathbb Z}_N)} \leq \frac{2 \epsilon}{1-2 \delta} \cdot \frac{1}{N} {||\widehat{f}||}_{L^1({\mathbb Z}_N)},\end{equation} 

and, consequently, 

\begin{equation} \label{equation:minerrorboundreversed_space} \frac{1}{|M|} \sum_{x \in M} |h(x)| \leq \frac{2 \varepsilon}{1 - 2\delta} \cdot \frac{1}{N} \sum_{x \in {\mathbb Z}_N} |f(x)|.\end{equation}

\end{theorem} 

Our next result is a more applicable version of Theorem \ref{theorem:quantitativeLogan}. Datasets/signals are generally noisy and setting $u(x)=f(x)$ for $x \notin M$ does not address this problem. To alleviate this issue, we alow a suitable $L^1$ difference between $u$ and $f$. 

\begin{theorem}
\label{theorem:looseLogan}
Let $f: {\mathbb Z}_N \to {\mathbb C}$, and suppose $f$ is transmitted but the values ${\{f(x)\}}_{x \in M}$ are unobserved. Suppose that $\widehat{f}$ is $L^1$-concentrated on $S \subset {\mathbb Z}_N$ with norm $\leq \epsilon$. Let $\delta=\frac{|M||S|}{N}$, $\delta' = \frac{|M^c| |S|}{N}$, and $\alpha \in [0, \frac{2 \varepsilon}{N \delta'}]$. Suppose $\delta < \frac{1}{2}$ and let
\begin{equation}
\label{equation:looseL1recoveryreversed}
g = argmin_u {||\widehat{u}||}_{L^1({\mathbb Z}_N)} \ \ \text{with the constraint} \ \| f - u \|_{L^1(M^c)} \leq \alpha \| f \|_{L^1(M^c)}.
\end{equation}

Let $h = f - g$. Then if $h \neq 0$, then

\begin{equation}
\label{equation:looseminerrorboundreversed}
\| \widehat{h} \|_{L^1(\mathbb{Z}_N)} \leq \frac{1}{N} \frac{2 \varepsilon + 2 N \alpha \delta'}{1 - 2\delta} \| \widehat{f} \|_{L^1(\mathbb{Z}_N)}
\end{equation}

and, consequently,

\begin{equation}
\label{equation:looseminerrorboundreversed_space}
\frac{1}{|M|} \| h \|_{L^1(\mathbb{Z}_N)} \leq \frac{2 \varepsilon + 2 N \alpha \delta'}{1 - 2\delta} \cdot \frac{1}{N} \| f \|_{L^1(\mathbb{Z}_N)}.
\end{equation}
\end{theorem}

\begin{remark}
It is worth noting that in equation (\ref{equation:looseL1recoveryreversed}), $f - u$, and therefore, $h$, can have support across all of $\mathbb{Z}_N$. 
\end{remark}

Theorems \ref{theorem:quantitativeLogan} and \ref{theorem:looseLogan} present results for any $M$. However, we can get much stronger, albeit probabilistic, results if we assume $M$ is generic in the following sense.

\begin{definition} \label{def:generic} Let $0< p <1$. Then, a random set $S \subset [n]= \{0, 1, \cdots, n-1\}$ is generic if each element of $[n]$ is selected independently with probability $p$.
\end{definition} 

We shall need the following result due to Talagrand (\cite{Talagrand98}),  which is a form of concentration inequality. 

\begin{theorem}[Talagrand's theorem] \label{thm: Talagrand98}  Let $(\varphi_j)_{j=1}^n$ be an orthonormal system in 
\( L^2(\Bbb Z_N)\) with \( \|\varphi_j\|_{L^\infty} \leq K \) for \( 1 \leq j \leq n \). There exists a constant $\gamma_0 \in (0,1)$ and a subset \( I \subset \{1, \dots, n\} \) with \( |I| \ge \gamma_0 n \) such that for every \( a = (a_i) \in \mathbb{C}^n \),

\[
\label{eq: Talagrand98}
\left( \sum_{i \in I} |a_i|^2 \right)^{1/2} \leq C_T K {( \log(n) \log \log(n))}^{\frac{1}{2}} {\left| \left| \sum_{i \in I} a_i \varphi_i \right| \right|}_{L_1},
\]

where \( C_T > 0 \) is a universal constant.
\end{theorem}

\begin{remark} \label{remark:talagrandexplanied} The insight from Talagrand's theorem suggests that even if the entire system has complex behavior, there is a large subset with ``good" concentration properties. 
Moreover, the proof shows that the conclusion holds with probability $1-o_N(1)$ for set $I$ satisfying the assumptions of Theorem \ref{thm: combo}, in the sense of Definition \ref{def:generic}. In our setting, we state Theorem \ref{thm: Talagrand98} as follows. \end{remark} 

\begin{theorem} \label{thm: combo} There exists $\gamma_0 \in (0,1)$ such that if $h: {\mathbb Z}_N \to {\mathbb C}$ supported in a generic set $M$ (in the sense of Definition \ref{def:generic}) of size $\gamma_0 \frac{N}{\log(N)}$, then with probability $1-o_N(1)$, 
\begin{equation} \label{eq: talagrand} {\left( \frac{1}{N}  \sum_{m \in {\mathbb Z}_N} {|\widehat{h}(m)|}^2 \right)}^{\frac{1}{2}} \leq C_T {(\log(N) \log \log(N))}^{\frac{1}{2}} \cdot \frac{1}{N} \sum_{m \in {\mathbb Z}_N} |\widehat{h}(m)|, \end{equation} where $C_T > 0$ is a constant that depends only on $\gamma_0$. 
\end{theorem} 

Using the notation in section \ref{subsection:notation}, equation \eqref{eq: talagrand} can be written as

$$\| \widehat{h} \|_{L^2(\mu)} \leq C_T {(\log(N) \log \log(N))}^{\frac{1}{2}} \| \widehat{h} \|_{L^1(\mu)}.$$

The proof of Theorem \ref{thm: combo} assumes \eqref{eq: talagrand} is true.  We apply Markov's inequality to show the equation is true over all $h$ supported in $M$ with probability at least $1- o_N(1)$.  Indeed, if each index of $I$ is selected with probability $\delta = \frac{1}{\log(N)}$, then a modification of the proof of Theorem \ref{thm: Talagrand98} on page 181 of  \cite{Talagrand98} gives us what we want.  Choosing $\delta$ in this way we have that,
$$
\mathbb{E}_{I}
\left (
\sup_{\|a\|_2 \leq 1} \left \| \sum_{i \in I} a_i\varphi_i\right\|
\right ) \leq \frac{K}{\sqrt{\log(1 / \delta)}}=
\frac{K}{\sqrt{\log \log N}},
$$
where $\| \cdot \|$ is the norm defined in on page 181 in \cite{Talagrand98} and the expectation is taken over $I$ defined in Definition \ref{def:generic}.  By Markov's inequality, we have that,
$$
\Pr \left (
\sup_{\|a\|_2 \leq 1} \left \| \sum_{i \in I} a_i\varphi_i\right\| \geq  \frac{1}{8} 
\right ) 
\leq \frac{K}{\sqrt{\log \log N}}.
$$
Hence with probability at least $1 - \frac{K}{\sqrt{\log \log N}} = 1 - o_N(1)$,
$$
\sup_{\|a\|_2 \leq 1} \left \| \sum_{i \in I} a_i\varphi_i\right\|
\leq \frac{1}{8},
$$
as needed.
\begin{remark} \label{rmk: talagrandconstant} The Talagrand constant, $C_T$, in Theorem \ref{thm: combo} depends on $\gamma_0$. The precise nature of this dependence is a very interesting area of exploration. 
\end{remark} 

Using Theorem \ref{thm: combo} we shall establish the following generic version of Theorem \ref{theorem:quantitativeLogan}. 

\begin{theorem} \label{theorem:quantitativerandomLogan} Let $f: {\mathbb Z}_N \to {\mathbb C}$, and suppose that the values ${\{f(x)\}}_{x \in M}$ are unobserved, where $M$ is a generic subset of ${\mathbb Z}_N$ (in the sense of Definition \ref{def:generic}), of size $\leq \gamma_0 \frac{N}{\log(N)}$, where $\gamma_0$ is as in Theorem \ref{thm: combo}. Let $g$ be as in (\ref{equation:simpleL1recoveryreversed}). Suppose that $\widehat{f}$ is $\epsilon$-concentrated on $S \subset {\mathbb Z}_N$ such that 
\begin{equation} \label{eq:sizespectral} |S|< \frac{1}{16 C_T^2} \frac{N}{\log(N) \log \log(N)}. \end{equation} 

Let $h=f-g$. Then if $h \not=0$, then with probability $1-o_N(1)$, 
\begin{equation} \label{equation:minerrorboundrandomreversed} {||\widehat{h}||}_{L^1({\mathbb Z}_N)} \leq \frac{4 \epsilon}{N} \cdot {||\widehat{f}||}_{L^1({\mathbb Z}_N)},\end{equation} 

and, consequently, 

\begin{equation} \label{equation:minerrorboundreversed_space2} \frac{1}{|M|} \sum_{x \in M} |h(x)| \leq 4 \epsilon \cdot \frac{1}{N} \sum_{x \in {\mathbb Z}_N} |f(x)|.\end{equation}

\end{theorem} 

\vskip.125in 

Much as with Theorem \ref{theorem:looseLogan}, we can take into account the fact that data is often noisy by allowing some difference between $g$ and $f$ on $M^c$.

\vskip.125in 

\begin{theorem} \label{theorem:quantitativerandomnoisyLogan} Let $f: {\mathbb Z}_N \to {\mathbb C}$, and suppose that the values ${\{f(x)\}}_{x \in M}$ are unobserved, where $M$ is a generic subset of ${\mathbb Z}_N$ (in the sense of Definition \ref{def:generic}), of size $\leq \gamma_0 \frac{N}{\log(N)}$, where $\gamma_0$ is as in Theorem \ref{thm: combo}. Let 
\begin{equation} \label{eq:noisyLogan} g = argmin_u {||\widehat{u}||}_1: {||u-f||}_{L^1(M^c)} \leq \delta N^{-1} {||f||}_{L^1(M^c)}. \end{equation}

Suppose that $\widehat{f}$ is $\epsilon$-concentrated on $S \subset {\mathbb Z}_N$ such that 
\begin{equation} \label{eq:sizespectralnoisy} |S|< \frac{1}{16 C_T^2} \frac{N}{\log(N) \log \log(N)}. \end{equation} 

Let $h=f-g$. Then if $h \not=0$, then with probability $1-o_N(1)$

\begin{equation} \label{equation:minerrorboundreversednoisy_space} \frac{1}{|M|} \sum_{x \in M} |h(x)| \leq \left( 4 \epsilon + 5 \delta \right) \cdot \frac{1}{N} \sum_{x \in {\mathbb Z}_N} |f(x)|.\end{equation}

\end{theorem} 

\vskip.125in 

\begin{remark}
Theorems \ref{theorem:quantitativerandomLogan} and \ref{theorem:quantitativerandomnoisyLogan} apply with high probability when the set of missing values, $M$, is generic with $\delta \leq \gamma_0$. The mechanism is explained following the statement of Theorem \ref{thm: combo}. 
\end{remark}

\vskip.125in

\begin{remark}
\label{remark:selecting_delta}
Theorems \ref{theorem:looseLogan} and \ref{theorem:quantitativerandomnoisyLogan} bound $\|g - f\|_{L^1(M^c)}$ by a parameter that needs to be specified. Experimentally, we've found good results by first taking
$$\tilde{f} = f - medfilt(f)$$
where $medfilt(\cdot)$ is a median filter and then requiring
$$\|g - f\|_{L^1(M^c)} \leq \frac{0.75}{|M^c|}\sum_{x \in M^c} |\tilde{f}(x) - m|$$
where $m$ is the median value of $\tilde{f}$ on $M^c$. However, more theoretical and experimental work is needed to improve this.
\end{remark}

\vskip.125in

\begin{remark} \label{remark:delta_eps_adding}
On their face, Theorems \ref{theorem:looseLogan} and \ref{theorem:quantitativerandomnoisyLogan} create an impression that defining $g$ according to equations \ref{equation:looseL1recoveryreversed} and \ref{eq:noisyLogan} would increase error. However, experimentation shows that the error generally decreases in this case. This makes sense as the ``loosening'' helps us avoid overfitting. More work is needed to understand this at a theoretical level.
\end{remark}

\vskip.125in

\begin{remark} \label{remark:almost_mae} Theorems \ref{theorem:quantitativeLogan} and \ref{theorem:quantitativerandomLogan} result in bounds on $\frac{1}{|M|} \sum_{x \in M} |h(x)|$. If it were the case that $\frac{1}{M} \sum_{x \in M} |f(x)| =\frac{1}{N} \sum_{x \in {\mathbb Z}_N} |f(x)|$, then we would be able to obtain bounds on

$$MAE_{\text{weighted}}(f,g) \equiv \frac{\sum_{x \in M} |h(x)|}{\sum_{x \in M} |f(x)|}.$$

Although we cannot obtain such an equality, we \textit{can} show that if $M$ is generic then

$$\left| \frac{1}{|M|} \sum_{x \in M} |f(x)| -\frac{1}{N} \sum_{x \in {\mathbb Z}_N} |f(x)| \right|$$

is very small with high probability. We turn our attention to this in Theorem \ref{theorem:hoeffding_M_vs_Z_N}.
\end{remark} 

\begin{remark}
Equation \ref{equation:minerrorboundreversednoisy_space} and Remark \ref{remark:almost_mae} indicate that that $MAE_{weighted}(f, g)$ is approximately bounded by $4\epsilon + 2\delta$. We might want, for example, $MAE_{weighted}(f, g) < 0.2$. If we make the simplifying assumption that $\frac{1}{|M|} \sum_{x \in M} |f(x)| = \frac{1}{N} \sum_{x \in \mathbb{Z}_N} |f(x)|$ (see Theorem \ref{theorem:hoeffding_M_vs_Z_N}) this puts relatively tight bounds on $\varepsilon$ and $\delta$ i.e.

$$\delta, \epsilon > 0 \,\,\,\,\, \text{and} \,\,\,\,\, 4\epsilon + 2\delta <0.2.$$

This restricts $\delta < 0.1$. Perhaps more consequentially, one can see that $f$ must be quite band-limited. For example, the above indicates that $\varepsilon < 0.05$ and so there must be a set, $S$, such that

$$\frac{|S|}{N} < \frac{1}{16 C_T^2} \frac{1}{\log(N) \log \log(N)} \overset{N \to \infty}{\longrightarrow} 0$$

but that, using Definition $\ref{def:eps_conc}$

$$\frac{\| \widehat{f} \|_{L^1(S)}}{\| \widehat{f} \|_{L^1(S^c)}} \geq \frac{N}{0.05} - 1 = 20N-1 \overset{N \to \infty}{\longrightarrow} 20N.$$

It is worth noting, however, just as in remark \ref{remark:delta_eps_adding}, experimentation shows that one can loosen these bounds and still obtain good results.
\end{remark}

\vskip.125in 

\begin{remark}
\label{remark:eps_concentrated}
Many of the above theorems take as an assumption that $\widehat{f}$ is $\varepsilon$-concentrated on some $S$ of sufficiently small size. Given a function $f$ with unobserved values, it is difficult to validate that $\widehat{f}$ satisfies this concentration condition. However, it is well known that large classes of reasonable signals are largely determined using a farily small number of Fourier coefficients. See, for example, Chapter 9 of \cite{Mallat2008-kc} where this matter is discussed in detail. 

We also note that the empirical results show that $L^1$-minimization often leads to accurate imputation even if the concentration conditions do not hold.
\end{remark}

\vskip.125in 

We now return to the point made in Remark \ref{remark:almost_mae}.

\vskip.125in 

\begin{theorem}
\label{theorem:hoeffding_M_vs_Z_N}
Assume $M$ is a generic subset of $\mathbb{Z}_N$ with probability $p$. Then, for large $N$,

\begin{equation}
\label{equation:hoeffding_M_vs_Z_N}
\Pr\Biggl(\Bigl|\frac{1}{|M|}\sum_{x\in M}|f(x)|-\frac{1}{N}\sum_{x\in\mathbb{Z}_N}|f(x)|\Bigr|\ge t\Biggr)
\le 2\left(1 - p + p\exp\left(-\frac{2 t^2}{\|f \|_{\infty}^2}\right)\right)^N.
\end{equation}
\end{theorem}

As the reader shall see, this result is proved using Hoeffding's inequality. Of particular importance is that, for any $t > 0$, the term inside the parentheses is in $[0, 1)$ and so

$$2 \left(1 - p + p\exp\left(-\frac{2 t^2}{\|f \|_{\infty}^2}\right)\right)^N \overset{N \to \infty}{\longrightarrow} 0.$$

\begin{remark}
Note that theorem \ref{theorem:hoeffding_M_vs_Z_N} is not very meaningful when $\| f \|_{\infty} \gg \| f\|_1$ e.g. when $f$ has aberrant spikes. As an extreme example, if $s_1 \leq |f(x)| \leq s_2$ for all $x \neq x_0 \in \mathbb{Z}_N$ and $f(x_0) \gg s_2$ then any generic $M$ is extremely unlikely to contain $x_0$ but and, as such,

$$\left| \frac{1}{|M|} \sum_{x \in M} |f(x)| - \frac{1}{N} \sum_{x \in \mathbb{Z}_N} |f(x)| \right|$$

is quite likely to be large.
\end{remark}

\begin{remark}
The theorems in this subsection can be easily generalized to functions from $\mathbb{Z}_N^d \to \mathbb{C}, d > 1$. The question of multivariate time-series is more nuanced and will be addressed in the sequel. 
\end{remark}

\subsubsection{Connections with the theory of exact signal recovery} 

In the theory of exact signal recovery, the basic question is whether the signal $f: {\mathbb Z}_N \to {\mathbb C}$ can be recovered from its Fourier transform if the frequencies ${\{\widehat{f}(m)\}}_{m \in S}$ are unobserved, where $S \subset {\mathbb Z}_N$. It was shown by Donoho and Stark (\cite{DS89}) and, independently, by Matolcsi and Szucs (\cite{MS73}) that if $E$ is the support of $f$, i.e $E=\{x \in {\mathbb Z}_N: f(x) \not=0\}$ and 

$$|E| \cdot |S|<\frac{N}{2},$$

then $f$ can be recovered exactly and uniquely. The roles of $f$ and $\widehat{f}$ can easily be interchanged, so if $f$ is transmitted directly, and the values ${\{f(x)\}}_{x \in M}$ are unobserved, then $f$ can be recovered uniquely and exactly if 

$$|M| \cdot |\{m \in {\mathbb Z}_N: \widehat{f}(m) \not=0 \}|<\frac{N}{2}.$$

A variety of recovery methods have been developed over the years. See, for example, \cite{DS89}, \cite{AS19}, \cite{KGR08}, \cite{KT07}, and the references contained therein. 

In the context of the problem considered in this paper, the set $M$ of missing values is random, so we can improve upon the recovery condition considerably using the following result from (\cite{IM24}) due to the second and the third listed authors using the celebrated Bourgain $\Lambda_q$ theorem (\cite{Bourgain89}). 

\begin{theorem} \label{signalrecoverycorollary}Let $f: {\mathbb Z}_N \to {\mathbb C}$ supported in $E \subset {\mathbb Z}_N$. Assume that $\{\hat f(m)\}_{m\in S}$ are unobserved,  where $S$ is a subset of ${\mathbb Z}_N^d$ of size $\lceil N^{\frac{2}{q}} \rceil$, for some $q>2$, randomly chosen with uniform probability. Then there exists a constant $C(q)$ that depends only on $q$ such that with probability $1-o_N(1)$, if 
\begin{equation} \label{recoverycondition} |E|< \frac{N}{2{(C(q)/\varepsilon)}^{\frac{1}{\frac{1}{2}-\frac{1}{q}}}}, \end{equation} then $f$ can be recovered uniquely where $\varepsilon = o_N(1)$. 
\end{theorem} 
The probabilistic part of Theorem \ref{signalrecoverycorollary} comes from Markov's inequality.  Indeed, from Bourgain's Theorem,
$$
\mathbb{E}_{I} 
\left (
\sup_{\|a\|_2 \leq 1} \left \|\sum_{i \in I} a_i \varphi_i  \right \|_q 
\right )
\leq C(q),
$$
where the randomness is on $I$ as in Definition \ref{def:generic}.  By Markov's inequality we have,
$$
\Pr \left (
\sup_{\|a\|_2 \leq 1} \left \|\sum_{i \in I} a_i \varphi_i \right \|_q  \geq \frac{C(q)}{\varepsilon}
\right ) \leq \frac{1}{\varepsilon},
$$
as needed.

\vskip.125in 

The recovery in Theorem \ref{signalrecoverycorollary} can be accomplished by the method least squares, i.e 
$$ f = arg min_u {||\widehat{u}-\widehat{f}||}_2: |spt(u)|=|spt(f)|,$$ following the procedure laid out by Donoho and Stark in \cite{DS89}. 

This is a very slow algorithm, so it is natural to check under what assumptions one can use Logan's algorithm (\ref{equation: simpleL1recoveryOrig}). This issue has already been thoroughly treated in Subsection \ref{subsection: L1 minimization and logan}, but we are now going to compare the threshold (\ref{recoverycondition}) with the corresponding condition needed to run Logan's algorithm. We begin with the following simple observation. In the case $q \ge 3$, it can be found in Vershynin's book. The general case follows from a similar argument, as was pointed out to us by William Hagerstrom (\cite{Hagerstrom2025}). 

\begin{lemma} \label{lemma:vershynintrick} Suppose that for $h: {\mathbb Z}_N^d \to {\mathbb C}$ with $\widehat{h}$ supported in $S \subset {\mathbb Z}_N^d$, 
\begin{equation} \label{eq:bourgainproxy} {\left( \frac{1}{N^d} \sum_{x \in {\mathbb Z}_N^d} {|h(x)|}^q \right)}^{\frac{1}{q}} \leq C(q) {\left( \frac{1}{N^d} \sum_{x \in {\mathbb Z}_N^d} {|h(x)|}^2 \right)}^{\frac{1}{2}} \end{equation} for some $q>2$. 

Then 
$$ {\left( \frac{1}{N^d} \sum_{x \in {\mathbb Z}_N^d} {|h(x)|}^2 \right)}^{\frac{1}{2}} \leq {(C(q))}^{\frac{q}{q-2}} \cdot \frac{1}{N^d} \sum_{x \in {\mathbb Z}_N^d} |h(x)|. $$
\end{lemma} 

To prove this result, observe that Holder's inequality implies that 
\begin{equation} \label{eq:vershynintrick}   {\left( \frac{1}{N^d} \sum_{x \in {\mathbb Z}_N^d} {|h(x)|}^2 \right)}^{\frac{1}{2}} \leq {\left( \frac{1}{N^d} \sum_{x \in {\mathbb Z}_N^d} |h(x)| \right)}^{\frac{1}{2}-\frac{1}{2(q-1)}} \cdot {\left( {\left( \frac{1}{N^d} \sum_{x \in {\mathbb Z}_N^d} {|h(x)|}^q \right)}^{\frac{1}{q}} \right)}^{\frac{q}{2(q-1)}}. \end{equation} 

Using (\ref{eq:bourgainproxy}) and (\ref{eq:vershynintrick}), we see that 
$$ {\left( \frac{1}{N^d} \sum_{x \in {\mathbb Z}_N^d} {|h(x)|}^2 \right)}^{\frac{1}{2}} \leq {\left( \frac{1}{N^d} \sum_{x \in {\mathbb Z}_N^d} |h(x)| \right)}^{\frac{1}{2}-\frac{1}{2(q-1)}} \cdot {(C(q))}^{\frac{q}{2(q-1)}} \cdot {\left( {\left( \frac{1}{N^d} \sum_{x \in {\mathbb Z}_N^d} {|h(x)|}^2 \right)}^{\frac{1}{2}} \right)}^{\frac{q}{2(q-1)}}$$ and the conclusion of Lemma \ref{lemma:vershynintrick} follows. 

\vskip.125in 

Using Lemma \ref{lemma:vershynintrick} we obtain the following result. 

\begin{theorem} \label{theorem:fromIMtoLogan} Let $S$ be a random subset of ${\mathbb Z}_N^d$ of size $\lceil N^{\frac{2d}{q}} \rceil$ for some $q>2$. Suppose that $f: {\mathbb Z}_N^d \to {\mathbb C}$ and the frequencies ${\{\widehat{f}(m)\}}_{m \in S}$ are unobserved. Suppose that $f$ is supported in $E \subset {\mathbb Z}_N^d$. If 
\begin{equation} \label{eq:transferencethreshold} |E|< \frac{N^d}{4 {(C(q)/\varepsilon)}^{\frac{2q}{q-2}}}, \end{equation} then $f$ can be recovered using Logan's method \eqref{equation: simpleL1recoveryOrig} with probability at least $1 - \varepsilon$. 
\end{theorem} 

Note that the recovery condition (\ref{eq:transferencethreshold}) is slightly weaker than the condition (\ref{recoverycondition}) above, but under the condition (\ref{eq:transferencethreshold}) we may use Logan's $L^1$ recovery method. 

\vskip.125in 

The role of $f$ and $\widehat{f}$ in Theorem \ref{signalrecoverycorollary} can once again be reversed to yield the following conclusion. 
\begin{theorem} \label{signalrecoverycorollaryreversed} Let $f: {\mathbb Z}_N \to {\mathbb C}$. Let $\Sigma=\left\{m \in {\mathbb Z}_N: \widehat{f}(m) \not=0 \right\}$. Assume that ${\{f(x)\}}_{x \in M}$ are unobserved, where $M$ is a subset of ${\mathbb Z}_N$ of size $\lceil N^{\frac{2}{q}} \rceil$, for some $q>2$, randomly chosen with uniform probability. Then there exists a constant $C(q)$ that depends only on $q$ such that with probability $1-o_N(1)$, if 
\begin{equation} \label{recoverycondition2} |\Sigma|< \frac{N}{2{(C(q)/\varepsilon)}^{\frac{1}{\frac{1}{2}-\frac{1}{q}}}}, \end{equation} then $\widehat{f}$ (and hence $f$) can be recovered uniquely.  Note $\varepsilon = o_N(1)$. 
\end{theorem} 

While the Fourier sparsity condition is not necessarily likely to hold, the Fourier transform of $f$ can be made sparse by force, introducing a small error, followed by an $L^1$ recovery procedure. See, for example, \cite{KT07}. We shall investigate this approach in detail in the sequel. 

\subsection{Numerical example}
\label{subsection: example}

In the companion paper, we will systematically test the algorithms described in Subsection \ref{subsection: L1 minimization and logan}. Here we briefly look at an example use of these algorithms. For this example, we've chosen the first $300$ records of a dataset which describes the monthly beer production in Australia\footnote{\url{https://www.kaggle.com/code/mpwolke/australian-monthly-beer-production}} and randomly removed half of the points. In this section, the ``linear interpolation'' algorithm refers to filling in the value, $f(x)$ with the value

$$f(a) + \frac{f(b) - f(a)}{b - a} x$$

where $b$ is the smallest value larger than $x$ not in $M$ and $a$ is the largest value smaller than $x$ not in $M$.

\begin{figure}
    \centering
    \includegraphics[width=0.49\linewidth]{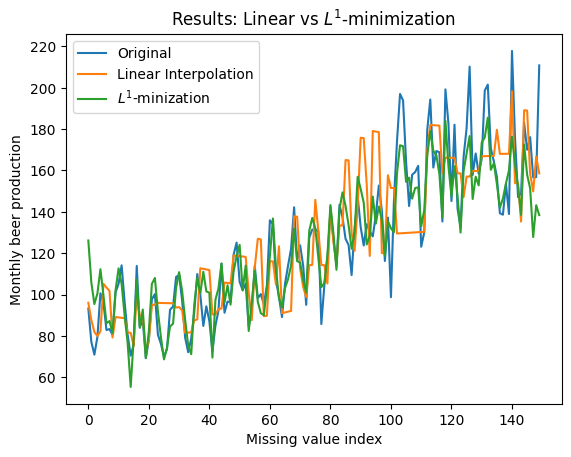}
    \includegraphics[width=0.49\linewidth]{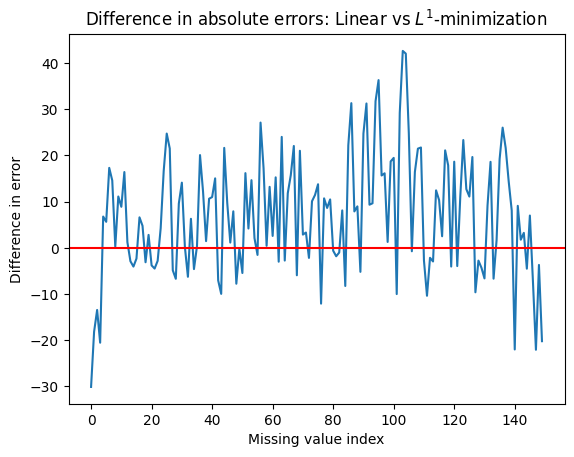}
    \caption{Linear vs $L^1$-minimization Results on Monthly Beer Production Data in Australia. The left figure shows the original missing values and the results of filling in the missing values with linear interpolation and $L^1$-minimization. The right figure shows difference in absolute errors with linear interpolation and $L^1$-minimization. Values above zero indicate that $L^1$-minimization improved results and values below zero indicate the opposite.}
    \label{fig:example_L1}
\end{figure}

To demonstrate the value of the $L^1$-minimization algorithm described in Theorem \ref{theorem:looseLogan}, we filled in the missing points using linear interpolation and the $L^1$-minimization algorithm. The results can be seen in Figure \ref{fig:example_L1}. The left figure shows the original missing values and the results of filling in the missing values with linear interpolation and $L^1$-minimization. The right figure shows difference in absolute errors between the two algorithms. In particular, the right plot shows

$$\{|f(x) - g_1(x)| - |f(x) - g_2(x)|\}_{x \in M}$$

where $f(x)$ is the true (original) data, $g_1(x)$ is the result of using linear interpolation to fill in the missing values, and $g_2(x)$ is the result of using $L^1$-minimization. Values above zero indicate that $L^1$-minimization improved results and values below zero indicate the opposite. 

The results can be most clearly seen in the right figure - the vast majority of the points are above zero, indicating that $L^1$-minimization improved the results on the vast majority of the points. Indeed, if we compare the $MAE$'s, we see that

$$\frac{MAE(g_2, f)}{MAE(g_1, f)} \approx 0.4171.$$

\begin{figure}
    \centering
    \includegraphics[width=0.49\linewidth]{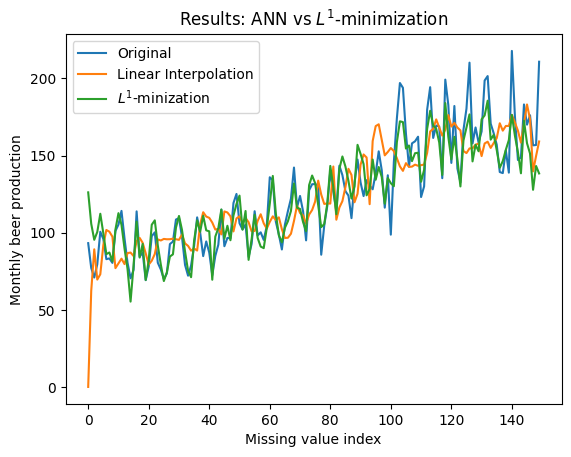}
    \includegraphics[width=0.49\linewidth]{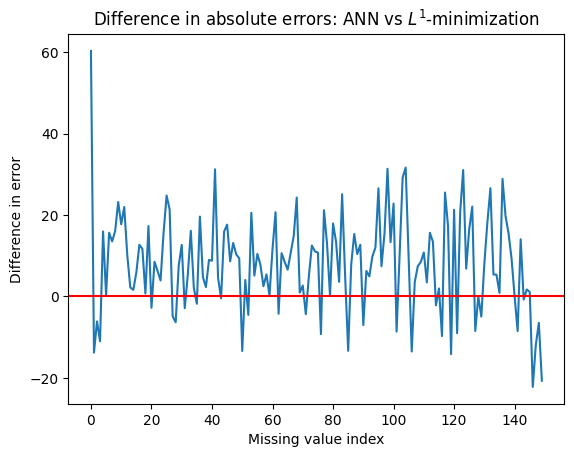}
    \caption{ANN vs $L^1$-minimization Results on Monthly Beer Production Data in Australia. The left figure shows the original missing values and the results of filling in the missing values with the ANN and $L^1$-minimization. The right figure shows difference in absolute errors with the ANN and $L^1$-minimization. Values above zero indicate that $L^1$-minimization improved results and values below zero indicate the opposite.}
    \label{fig:example_ANN_vs_L1}
\end{figure}

Similarly, we conduct the same test with $g_1$ as the result of using an artificial neural network (``ANN'') instead of linear interpolation.\footnote{This network has two hidden layers, each with $64$ neurons. Each training run trains with the non-missing values with the mean squared error loss function and the Adam (\cite{Kingma2014Adam}) optimizer for five epochs.} The results can be seen in Figure \ref{fig:example_ANN_vs_L1}. As can be seen, the $L^1$-optimizer performs better than the ANN on the vast majority of points. If we compare the MAE's we get

$$\frac{MAE(g_2, f)}{MAE(g_1, f)} \approx 0.4646$$

which is only slightly worse than the comparison with linear interpolation.

\vskip.125in 

\section{Proofs of theorems}
\label{section: proofs}

\subsection{Proof of Theorem \ref{theorem:quantitativeLogan}} We have 
$$ {||\widehat{g}||}_{L^1({\mathbb Z}_N)} = {||\widehat{f}-\widehat{h}||}_{L^1({\mathbb Z}_N)}$$
$$ = {||\widehat{f} - \widehat{h}||}_{L^1(S)} + {||\widehat{h} - \widehat{f}||}_{L^1(S^c)} \ge {|| \widehat{f} ||}_{L^1(S)} - {|| \widehat{h} ||}_{L^1(S)} + {|| \widehat{h} ||}_{L^1(S^c)}-{|| \widehat{f} ||}_{L^1(S^c)}. $$

Since $\widehat{f}$ is $\epsilon$-concentrated on $S$, the right hand side is bounded from below by 

$${|| \widehat{f} ||}_{L^1(S)}- \frac{\epsilon}{N} {|| \widehat{f} ||}_{L^1({\mathbb Z}_N)} + \left( || \widehat{h} ||_{L^1(S^c)} - || \widehat{h}||_{L^1(S)} \right)$$
\begin{equation} \label{equation:almost} \ge \left( 1 - 2 \frac{\epsilon}{N} \right) || \widehat{f} ||_{L^1({\mathbb Z}_N)} + \left(|| \widehat{h} ||_{L^1(S^c)} - || \widehat{h} ||_{L^1(S)}\right).\end{equation} 

We now observe that 
$$ | \widehat{h}(m)| = \left| N^{-\frac{1}{2}} \sum_{x \in M} \chi(-x \cdot m) {h}(x) \right|\leq N^{-1} |M| \cdot {|| \widehat{h} ||}_{L^1({\mathbb Z}_N)}.$$

It follows that 
$$ {|| \widehat{h} ||}_{L^1(S)} \leq \delta \cdot {|| \widehat{h}||}_{L^1({\mathbb Z}_N)},$$ 

hence

$${|| \widehat{h} ||}_{L^1(S^c)}-{|| \widehat{h} ||}_{L^1(S)}={|| \widehat{h} ||}_{L^1({\mathbb Z}_N)} - 2{|| \widehat{h} ||}_{L^1(E)}$$
$$ \ge (1-2 \delta) {|| \widehat{h} ||}_{L^1({\mathbb Z}_N)}.$$

Plugging this back into (\ref{equation:almost}), we see that 

$${|| \widehat{g}||}_{L^1({\mathbb Z}_N)} \ge {|| \widehat{f} ||}_{L^1({\mathbb Z}_N)}+(1-2 \delta) {|| \widehat{h} ||}_{L^1({\mathbb Z}_N)}-2 \frac{\epsilon}{N} {|| \widehat{f} ||}_{L^1({\mathbb Z}_N)}.$$

To avoid contradiction the fact that $g$ is the $L^1$ minimizer, we must have 

$$(1-2 \delta) {|| \widehat{h} ||}_{L^1({\mathbb Z}_N)} - 2 \frac{\epsilon}{N} {||\widehat{f}||}_{L^1({\mathbb Z}_N)} \leq 0,$$

so 

$${||\widehat{h}||}_{L^1({\mathbb Z}_N)} \leq \frac{2 \epsilon}{1 - 2 \delta} \cdot \frac{1}{N} {|| \widehat{f} ||}_{L^1({\mathbb Z}_N)}.$$ 

To obtain equation (\ref{equation:minerrorboundreversed_space}) we note that for any function $\alpha:\mathbb{Z}_N \to \mathbb{C}$ we know that for $x \in {\mathbb Z}_N,$

$$\alpha(x) = N^{-\frac{1}{2}} \sum_{m \in \mathbb{Z}_N^d } \chi(xm) \widehat{\alpha} (m).$$

So,

$$|\alpha(x)| \leq N^{-\frac{1}{2}}  \| \widehat{\alpha} \|_{L^1({\mathbb Z}_N)} \Rightarrow \| \alpha \|_{L^1({\mathbb Z}_N)} \leq |supp(\alpha)| \cdot  N^{-\frac{1}{2}} \| \widehat{\alpha} \|_{L^1({\mathbb Z}_N)}.$$

Similarly, 

$$\|\widehat{\alpha}\|_{L^1({\mathbb Z}_N)} \leq |supp(\widehat{\alpha})| \cdot  N^{-\frac{1}{2}}  \| \alpha \|_{L^1({\mathbb Z}_N)}.$$

So, we get

$$\frac{N^{\frac{1}{2}}}{|M|} \|h\|_{L^1({\mathbb Z}_N)} \leq \frac{2 \varepsilon}{1 - 2\delta} \frac{N^{\frac{1}{2}}}{N} \|f\|_{L^1({\mathbb Z}_N)}
\Rightarrow \frac{1}{|M|} \|h\|_{L^1({\mathbb Z}_N)} \leq \frac{2 \varepsilon}{1 - 2\delta} \cdot \frac{1}{N} \|f\|_{L^1({\mathbb Z}_N)}.$$

\subsection{Proof of Theorem \ref{theorem:looseLogan}}
Just as in the proof of theorem \ref{theorem:quantitativeLogan}, we have equation (\ref{equation:almost}).  Now, know that for $m \in \mathbb{Z}_N$,
$$\widehat{h}(m) = N^{-\frac{1}{2}} \sum_{x \in M} \chi(x \cdot m) h(x) + N^{-\frac{1}{2}} \sum_{x \notin M} \chi(x \cdot m) h(x)$$
and so
\begin{align*}
| \widehat{h} (m) |
& \leq \frac{|M|}{N} \| \widehat{h} \|_{L^1(\mathbb{Z}_N)} + \frac{1}{N^{\frac{1}{2}}} \| h \|_{L^1(M^c)} \\
& \leq \frac{|M|}{N} \| \widehat{h} \|_{L^1(\mathbb{Z}_N)} + \frac{\alpha}{N^{\frac{1}{2}}} \| f \|_{L^1(M^c)} \\
& \leq \frac{|M|}{N} \| \widehat{h} \|_{L^1(\mathbb{Z}_N)} + \frac{\alpha |M^c|}{N} \| \widehat{f} \|_{L^1(\mathbb{Z}_N)}.
\end{align*}
Thus,
$$ \| \widehat{h} \|_{L^1(S)} \leq \delta \| \widehat{h} \|_{L^1(\mathbb{Z}_N)} + \alpha \delta' \| \widehat{f} \|_{L^1(\mathbb{Z}_N)}.$$

Plugging this back into equation (\ref{equation:almost}) gets

$$\| \widehat{g} \|_{L^1(\mathbb{Z}_N)} \geq \| \widehat{f} \|_{L^1(\mathbb{Z}_N)} + \| \widehat{h} \|_{L^1(\mathbb{Z}_N)} (1 - 2 \delta) + \| \widehat{f} \|_{L^1(\mathbb{Z}_N)} \left( -2 \alpha \delta ' - 2 \frac{\varepsilon}{N} \right).$$

Again, to avoid a contradiction, we need

$$\| \widehat{f} \|_{L^1(\mathbb{Z}_N)} + \| \widehat{h} \|_{L^1(\mathbb{Z}_N)} (1 - 2 \delta) + \| \widehat{f} \|_{L^1(\mathbb{Z}_N)} \left( -2 \alpha \delta ' - 2 \frac{\varepsilon}{N} \right) \leq \| \widehat{f} \|_{L^1(\mathbb{Z}_N)}$$ 
$$\Leftrightarrow \| \widehat{h} \|_{L^1(\mathbb{Z}_N)} \leq \frac{2\frac{\varepsilon}{N} + 2 \alpha \delta'}{1 - 2 \delta} \| \widehat{f} \|_{L^1(\mathbb{Z}_N)} = \frac{1}{N}\frac{2 \varepsilon + 2 N \alpha \delta'}{1 - 2 \delta} \| \widehat{f} \|_{L^1(\mathbb{Z}_N)}.$$

Again, as in the proof of theorem \ref{theorem:quantitativeLogan}, we get

$$\frac{1}{|M|} \| h \|_{L^1(\mathbb{Z}_N)} \leq \frac{2 \varepsilon + 2 N \alpha \delta'}{1 - 2\delta} \cdot \frac{1}{N} \| f \|_{L^1(\mathbb{Z}_N)}.$$

\subsection{Proof of Theorem \ref{theorem:quantitativerandomLogan}} As usual, let $g$ denote the quantity in \ref{equation:simpleL1recoveryreversed}. Let $f=g+h$. We have 
$$ {||\widehat{g}||}_1 = {||\widehat{f}-\widehat{h}||}_1$$
$$={||\widehat{f}-\widehat{h}||}_{L^1(S)}+{||\widehat{f}-\widehat{h}||}_{L^1(S^c)}$$
$$ \ge \left( {||\widehat{f}||}_{L^1(S)}-{||\widehat{f}||}_{L^1(S^c)} \right) + \left( {||\widehat{h}||}_{L^1(S^c)} -{||\widehat{h}||}_{L^1(S)} \right)$$
$$ \ge {||\widehat{f}||}_1 \cdot \left(1-\frac{2\epsilon}{N} \right) + \left( {||\widehat{h}||}_{L^1(S^c)} -{||\widehat{h}||}_{L^1(S)} \right),$$ where we have used the concentration assumption in the last line. 

There are two possibilities. If 
$$ -\frac{2\epsilon}{N} {||\widehat{f}||}_1  + \left( {||\widehat{h}||}_{L^1(S^c)} -{||\widehat{h}||}_{L^1(S)} \right)>0,$$ we contradict the fact that $\widehat{g}$ is the $L^1$ minimizer, which forces $h$ to be identically $0$. Otherwise, we must have 
\begin{equation} \label{eq: dnieper}  {||\widehat{h}||}_{L^1(S^c)} -{||\widehat{h}||}_{L^1(S)} \leq \frac{2\epsilon}{N} {||\widehat{f}||}_1. \end{equation} 

To exploit this, observe that 
$$ {||\widehat{h}||}_{L^1(S)} \leq {|S|}^{\frac{1}{2}} \cdot {||\widehat{h}||}_{L^2(S)} \leq {|S|}^{\frac{1}{2}} \cdot {||\widehat{h}||}_2 $$
$$ = {|S|}^{\frac{1}{2}} \cdot N^{\frac{1}{2}} \cdot {||\widehat{h}||}_{L^2(\mu)}$$
$$ \leq {|S|}^{\frac{1}{2}} \cdot N^{-\frac{1}{2}} \cdot C_T \sqrt{\log(N) \log \log(N)} \cdot {||\widehat{h}||}_1,$$
where in the last line we applied Talagrand's bound (Theorem \ref{thm: combo} above). Plugging this back into (\ref{eq: dnieper}), we obtain 
$$ {||\widehat{h}||}_1 (1-2{|S|}^{\frac{1}{2}} \cdot N^{-\frac{1}{2}} \cdot C_T \sqrt{\log(N) \log \log(N)}) \leq \frac{2 \epsilon}{N} {||\widehat{f}||}_1$$ with probability $1-o_N(1)$. 

Using the assumption regarding $|S|$, this implies that 

$$ {||\widehat{h}||}_1 \leq \frac{4 \epsilon}{N} {||\widehat{f}||}_1, $$ 

which is (\ref{equation:minerrorboundrandomreversed}), once again, with probability $1-o_N(1)$.  To prove (\ref{equation:minerrorboundreversed_space2}), we again proceed as we did in the proofs of theorems \ref{theorem:quantitativeLogan} and \ref{theorem:looseLogan} and get

$$\frac{1}{|M|} \| h \|_{L^1(\mathbb{Z}_N)} \leq 4\epsilon \cdot \frac{1}{N} \| f \|_{L^1(\mathbb{Z}_N)}.$$

This completes the proof of Theorem \ref{theorem:quantitativerandomLogan}. 

\vskip.25in 

\subsection{Proof of Theorem \ref{theorem:quantitativerandomnoisyLogan}} 

\vskip.125in 

Let $f=g+h$. We have 
$$ {||\widehat{g}||}_1 = {||\widehat{f}-\widehat{h}||}_1$$
$$={||\widehat{f}-\widehat{h}||}_{L^1(S)}+{||\widehat{f}-\widehat{h}||}_{L^1(S^c)}$$
$$ \ge \left( {||\widehat{f}||}_{L^1(S)}-{||\widehat{f}||}_{L^1(S^c)} \right) + \left( {||\widehat{h}||}_{L^1(S^c)} -{||\widehat{h}||}_{L^1(S)} \right)$$
$$ \ge {||\widehat{f}||}_1 \cdot \left(1-\frac{2\epsilon}{N} \right) + \left( {||\widehat{h}||}_{L^1(S^c)} -{||\widehat{h}||}_{L^1(S)} \right),$$ where we have used the concentration assumption in the last line. 

There are two possibilities. If 
$$ -\frac{2\epsilon}{N} {||\widehat{f}||}_1  + \left( {||\widehat{h}||}_{L^1(S^c)} -{||\widehat{h}||}_{L^1(S)} \right)>0,$$ we contradict the fact that $\widehat{g}$ is the $L^1$ minimizer, which forces $h$ to be identically $0$. Otherwise, we must have 
\begin{equation} \label{eq: dniepernoisy} {||\widehat{h}||}_{L^1(S^c)} -{||\widehat{h}||}_{L^1(S)} \leq \frac{2\epsilon}{N} {||\widehat{f}||}_1. \end{equation} 

To exploit this, we write 
$$ h(x)=1_Mh(x)+1_{M^c}h(x). $$

This is necessary because $h$ is not supported in $M$ as it was in the previous theorem. We first deal with ${||\widehat{1_{M^c}h}||}_{L^1(S)}$. This quantity is bounded by 
$$ |S| \cdot {||\widehat{1_{M^c}h}||}_{\infty} \leq |S| \cdot N^{-\frac{1}{2}} \cdot {||1_{M^c}h||}_1.$$ 

By assumption (\ref{eq:noisyLogan}), this quantity is bounded by
$$ |S| \cdot N^{-\frac{1}{2}} \cdot \delta \cdot N^{-1} {||f||}_{L^1(M^c)} \leq |S| \cdot N^{-\frac{1}{2}} \cdot \delta \cdot {||f||}_{L^1(\mu)}.$$

We now deal with ${||\widehat{1_Mh}||}_{L^1(S)}$. We have 
$$ {||\widehat{1_Mh}||}_{L^1(S)} \leq {|S|}^{\frac{1}{2}} \cdot {||\widehat{1_Mh}||}_{L^2(S)} \leq {|S|}^{\frac{1}{2}} \cdot {||\widehat{1_Mh}||}_2 $$
$$ = {|S|}^{\frac{1}{2}} \cdot N^{\frac{1}{2}} \cdot {||\widehat{1_Mh}||}_{L^2(\mu)}$$
$$ \leq {|S|}^{\frac{1}{2}} \cdot N^{-\frac{1}{2}} \cdot C_T \sqrt{\log(N) \log \log(N)} \cdot {||\widehat{1_Mh}||}_1$$ 
$$\leq \frac{{||\widehat{1_Mh}||}_1}{4} \leq \frac{{||\widehat{1_{M^c}h}||}_1 + \|\widehat{h}\|_1}{4}$$
$$\leq \frac{N^{\frac{1}{2}} \cdot \delta \cdot {||f||}_{L^1(\mu)} + \|\widehat{h}\|_1}{4}$$
where in the third to last line we applied Talagrand's bound (Theorem \ref{thm: combo} above) and in the second to last line we use the assumption about $|S|$. 

It follows that the left-hand side of (\ref{eq: dniepernoisy}) is bounded from below by 
$$\|\widehat{h}\|_1 - \frac{5}{2} \cdot N^{\frac{1}{2}} \cdot \delta \cdot {||f||}_{L^1(\mu)} - \frac{\|\widehat{h}\|_1}{2} = \frac{\|\widehat{h}\|_1}{2} - \frac{5}{2}\cdot N^{\frac{1}{2}} \cdot \delta \cdot {||f||}_{L^1(\mu)}$$

in view of (\ref{eq:sizespectralnoisy}). 

Plugging everything back into (\ref{eq: dniepernoisy}), we see that 
\begin{equation} \label{eq: almostnitsahon} \|\widehat{h}\|_1 \leq \frac{4\epsilon}{N} {||\widehat{f}||}_1 + 5 \cdot N^{\frac{1}{2}} \cdot \delta \cdot {||f||}_{L^1(\mu)}\end{equation} 

\vskip.125in 

To prove (\ref{equation:minerrorboundreversednoisy_space}), we write 
$$ \frac{1}{M} \sum_{x \in M} |h(x)| \leq N^{-\frac{1}{2}} {||\widehat{h}||}_1.$$

Applying (\ref{eq: almostnitsahon}), we see that this quantity is bounded by

$$\left( 4\epsilon + 5 N \cdot N^{-1} \cdot \delta \right) \frac{1}{N} \sum_x |f(x)| \leq \left( 4\epsilon + 5 \delta \right) \frac{1}{N} \sum_x |f(x)|.$$


\subsection{Proof of Theorem \ref{theorem:hoeffding_M_vs_Z_N}}

We shall need the following concentration bound, known as Hoeffding's Inequality. 

\begin{theorem}[Hoeffding \cite{Hoeffding63}] \label{thm:hoeffding_finite_population} Let $X = \{x_1, ..., x_N\}$ be a finite population of real numbers and, for some $0 \leq n \leq N$, let $M \subset X$ be a random subset of $n$ elements where the elements are chosen without replacement and $S = n^{-1}\sum_{y \in M} y$. Then, for all $t > 0$

\[
\mathbb{P} \left( S - \mu \geq t \right) \leq \exp \left( -\frac{2nt^2}{(b - a)^2} \right).
\]
Furthermore,
\[
\mathbb{P} \left( |S - \mu| \geq t \right) \leq 2\exp \left( -\frac{2nt^2}{(b - a)^2} \right).
\]

Here, $\mu = n^{-1} \sum_{x \in X} x$, $b = \max X$ and $a = \min X$.

\end{theorem}

We can use theorem \ref{thm:hoeffding_finite_population} below. Let $A \sim Bin(N, p)$.

\begin{align*}
& Pr\left( \left| \frac{1}{|M|} \sum_{x \in M} |f(x)| - \frac{1}{N} \sum_{x \in \mathbb{Z}_N} |f(x)| \right| \geq t\right) \\
& = \sum_{k = 0}^{N} Pr \left( |M| = k\right) \cdot  Pr\left( \left| \frac{1}{|M|} \sum_{x \in M} |f(x)| - \frac{1}{N} \sum_{x \in \mathbb{Z}_N} |f(x)| \right| \geq t : |M| = k\right) \\
& = \sum_{k = 0}^{N} {N \choose k} p^k (1 - p)^{N - k} \cdot  Pr\left( \left| \frac{1}{|M|} \sum_{x \in M} |f(x)| - \frac{1}{N} \sum_{x \in \mathbb{Z}} |f(x)| \right| \geq t : |M| = k\right) \\
& \leq \sum_{k = 0}^{N} {N \choose k} p^k (1 - p)^{N - k} \cdot 2\exp \left( -\frac{2kt^2}{\|f \|_{\infty}^2} \right) \\
& = 2 \mathbb{E}\left[\exp \left( -\frac{2 A t^2}{\|f \|_{\infty}^2} \right)\right]
\end{align*}

where we used Hoeffding's inequality to obtain the fourth line. The expectation above is just the moment generating function for the binomial distribution with parameter $-\frac{2 t^2}{\|f \|_{\infty}^2}$ which is

$$\left(1 - p + p\exp\left({-\frac{2 t^2}{\|f \|_{\infty}^2}}\right)\right)^N.$$




%

\vskip.250in

\section{Appendix: Proof of Theorem \ref{theorem:DSLogan89Orig}}
\label{section:appendix}

In this section, we provide the proof of Theorem $\ref{theorem:DSLogan89Orig}$, originally found in \cite{DS89}, for completeness. 

\begin{proof}
Let $g$ be the solution of (\ref{equation: simpleL1recoveryOrig}). Then $f=g+h$, where $\widehat{h}$ is supported in $S$. We have 

$${||g||}_{L^1({\mathbb Z}_N)} = {||f-h||}_{L^1({\mathbb Z}_N)}$$
$$ = {||f-h||}_{L^1(E)}+{||h||}_{L^1(E^c)} \ge {||f||}_{L^1({\mathbb Z}_N)}+\left({||h||}_{L^1(E^c)}-{||h||}_{L^1(E)} \right).$$ 

Assume $g \not\equiv f$. If we can show that ${||h||}_{L^1(E^c)}>{||h||}_{L^1(E)}$ we will have shown that ${||g||}_{L^1({\mathbb Z}_N)}>{||f||}_{L^1({\mathbb Z}_N)}$, and obtained a contradiction since $g$ has the minimal $L^1$ norm. To obtain this inequality for $h$, observe that 
$$|h(x)|= \left| N^{-\frac{1}{2}} \sum_{m \in S} \chi(x \cdot m) \widehat{h}(m) \right|$$ 
$$= \left| N^{-\frac{1}{2}} \sum_{m \in S} \chi(x \cdot m) \cdot N^{-\frac{1}{2}}\sum_{y \in \mathbb{Z}_N} \chi(-y \cdot m) h(y) \right| \leq N^{-1} \cdot |S| \cdot {||h||}_{L^1({\mathbb Z}_N)}.$$
It follows that 
$${||h||}_{L^1(E)} \leq \frac{|E||S|}{N} \cdot {||h||}_{L^1({\mathbb Z}_N)}.$$

Since $|E||S|<\frac{N}{2}$ by assumption, it follows that ${||h||}_{L^1(E^c)} >{||h||}_{L^1(E)}$, as desired. 
\end{proof}

\vskip.250in

\section{Statements and Declarations}

\subsection{Competing Interests}
On behalf of all authors, the corresponding author states that there is no conflict of interest.

\subsection{Data Availability}
The dataset on monthly beer production in Australia is available at \url{https://www.kaggle.com/code/mpwolke/australian-monthly-beer-production}

\subsection{Author Contributions}
W. Burstein, A. Iosevich, A. Mayeli, and H. Nathan contributed equally to this work.


\newpage 

\begin{thebibliography}{8}

\bibitem{AS19} L. Abreu and M. Speckbacher, {\it Deterministic guarantees for $L^1$-reconstruction: a large sieve approach with geometric flexibility}, 13th International conference on Sampling Theory and Applications (SampTA), (2019).


\bibitem{Bourgain89} J. Bourgain, {\it Bounded orthogonal systems and the $\Lambda(p)$-set problem}, Acta Math. \textbf{162} (1989), no. 3-4, 227–245.


\bibitem{DS89} D. Donoho and P. Stark, {\it Uncertainty principle and signal processing}, SIAM Journal of Applied Math., (1989), Society for Industrial and Applied Mathematics, volume 49, No. 3, pp. 906-931. 



\bibitem{Hagerstrom2025} W. Hagerstrom, {\it A number of perspectives on signal recovery}, University of Rochester Honors Thesis (2025). 

\bibitem{Hoeffding63} W. Hoeffding, {\it Probability inequalities for sums of bounded random variables}, Journal of the American Statistical Association \textbf{58} (301): 13–30.


\bibitem{IKLM24} A. Iosevich, B. Kashin, I. Limonova, and A. Mayeli, {\it Subsystems of orthogonal systems and the recovery of sparse signals in the presence of random losses}, Russian Mathematical Surveys, (2024), Volume 79, Issue 6, Pages 1095-1097.

\bibitem{IM24}
A. Iosevich and A. Mayeli, {\it Uncertainty Principles, Restriction, Bourgain's $\Lambda_q$ theorem, and Signal Recovery},  Applied and Computational Harmonic Analysis 76 (2025): 101734.



\bibitem{KT07} B. Kashin and V. Temlyakov, {\it A remark on the problem of compressed sensing}, (Russian) Mat. Zametki 82 (2007), no. 6, 829–837; translation in Math. Notes 82 (2007), no. 5-6, 748-755.

\bibitem{Kingma2014Adam} D. Kingma and J. Ba, {\it Adam: A Method for Stochastic Optimization}, CoRR, abs/1412.6980, (2014).

\bibitem{KGR08} F. Krahmer, E. Gotz, and P. Rashkov, {\it Uncertainty in time-frequency representations on finite abelian groups and applications}, Appl. Comput. Harmon. Anal. \textbf{25} (2008), no. 2, 209-225.


\bibitem{Mallat2008-kc}
S.~Mallat, \emph{A Wavelet Tour of Signal Processing}, 3rd ed.,  
edited by G.~Peyré, Academic Press, Dec. 2008.

\bibitem{MS73} T. Matolcsi and J. Szucs, {\it Intersection des mesures spectrales conjug\'ees.} C.R. Acad. Sci. S\'er. I Math. \textbf{277} (1973), 841-843.

\bibitem{Smith07} J. Smith, {\it Mathematics of the Discrete Fourier Transform (DFT): with Audio Applications}, Second Edition, Stanford, (2007). 


\bibitem{Talagrand98} M. Talagrand, {\it Selecting a proportion of characters}, Israel J. Math. \textbf{108} (1998),
173-191.



\end{thebibliography}
\end{document}